\documentclass[11pt,twoside]{amsart}
\usepackage{latexsym,amssymb,amsmath}

\numberwithin{equation}{section}
\newtheorem{theorem}{Theorem}[section]
\newtheorem{lemma}[theorem]{Lemma}
\newtheorem{proposition}[theorem]{Proposition}
\newtheorem{corollary}[theorem]{Corollary}
\newtheorem{conjecture}[theorem]{Conjecture}

\theoremstyle{definition}
\newtheorem{definition}[theorem]{Definition} 
 
\newtheorem{remark}[theorem]{Remark}
\newtheorem{example}[theorem]{Example}

\begin{document}


\newcommand{\m}[1]{\marginpar{\addtolength{\baselineskip}{-3pt}{\footnotesize
\it #1}}}
\newcommand{\A}{\mathcal{A}}
\newcommand{\K}{\mathcal{K}}
\newcommand{\CC}{\mathcal{C}}
\newcommand{\knd}{\mathcal{K}^{[d]}_n}
\newcommand{\F}{\mathcal{F}}
\newcommand{\N}{\mathbb{N}}
\newcommand{\pr}{\mathbb{P}}
\newcommand{\Z}{\mathbb{Z}}
\newcommand{\R}{\mathbb{R}}
\newcommand{\I}{\mathit{I}}
\newcommand{\G}{\mathcal{G}}
\newcommand{\D}{\mathcal{D}}
\newcommand{\x}{\mathbf{x}}
\newcommand{\lcm}{\operatorname{lcm}}
\newcommand{\ndp}{N_{d,p}}
\newcommand{\tor}{\operatorname{Tor}}
\newcommand{\reg}{\operatorname{reg}}
\newcommand{\mf}{\mathfrak{m}}

\def\bb{{{\rm \bf b}}}
\def\cc{{{\rm \bf c}}}

\title[Ehrhart Clutters]{Ehrhart clutters: Regularity and Max-Flow Min-Cut}

\author{Jos\'e Mart\'{\i}nez-Bernal, Edwin O'Shea and Rafael H. Villarreal}
\footnote{The second author was partially supported by SNI. The third
author was partially supported
by CONACyT grant 49251-F and SNI.} 
\address{
Departamento de
Matem\'aticas\\
Centro de Investigaci\'on y de Estudios Avanzados del
IPN\\
Apartado Postal
14--740 \\
07000 Mexico City, D.F. } \email{[jmb,edwin,vila]@math.cinvestav.mx}


\begin{abstract}
If $\CC$ is a clutter with $n$ vertices and $q$ edges whose 
clutter matrix has column vectors ${\mathcal A} = \{v_1, \ldots, v_q\}$, 
we call $\CC$ an {\it Ehrhart clutter\/} if 
$\{(v_1,1),\ldots,(v_q,1)\} \subset \{ 0,1 \}^{n+1}$ 
is a Hilbert basis. 
Letting $A(P)$ be the Ehrhart ring of 
$P={\rm conv}(\mathcal{A})$, we are able to show that if
$\mathcal{C}$ is a uniform unmixed MFMC 
clutter, then $\CC$ is an Ehrhart clutter and in this 
case we provide sharp upper bounds on the Castelnuovo-Mumford regularity 
and the $a$-invariant of $A(P)$. Motivated by the Conforti-Cornu{\' e}jols conjecture on
packing problems, 
we conjecture that if $\CC$ is both ideal and the clique clutter 
of a perfect graph, then $\CC$ has the MFMC property. We prove this 
conjecture for Meyniel graphs, by showing that the clique clutters 
of Meyniel graphs are Ehrhart clutters. In much the same spirit, we 
provide a simple proof of our conjecture when $\CC$ is a uniform clique 
clutter of a perfect graph. We close with a generalization 
of Ehrhart clutters as it relates to total dual integrality.
\end{abstract}

\maketitle
\section{Introduction} 

A {\it clutter\/} $\CC$ is a family $E$ of subsets of a
finite ground set $X$ such that if $S_1, S_2 \in E$, then $S_1\not\subset S_2$. 
The ground set $X$ is called the {\em vertex set} of $\CC$ and $E$ 
is called the {\em edge set} of $\CC$, they are denoted by $V(\CC)$
and $E(\CC)$  respectively. Clutters are special hypergraphs and are sometimes called 
{\em Sperner families} in the literature. We can also think of a clutter 
as the maximal faces of a simplicial complex over a ground set. One 
example of a clutter is a graph with the vertices and edges defined in the 
usual way for graphs. For a thorough study of clutters and
hypergraphs from the point of view of combinatorial optimization 
and commutative algebra see \cite{cornu-book,Schr2}
and  \cite{faridijct,reesclu,tai-vantuyl} respectively.

Let $\mathcal{C}$ be a clutter with vertex set $X=\{x_1,\ldots,x_n\}$ and 
with edge set $E(\CC)$. 
We shall assume that $\mathcal{C}$ has no isolated vertices, 
i.e., each vertex occurs in at least one edge and every edge contains 
at least two vertices. 
Permitting an abuse of notation, we will also denote by $x_i$ the $i^{\textup{th}}$ 
variable in the polynomial ring $R=K[x_1,\ldots,x_n]$ over a field $K$. 
The {\it edge ideal\/} of $\mathcal{C}$, denoted by $I(\mathcal{C})$, is the 
monomial ideal of $R$ generated by all 
monomials $x_e=\prod_{x_i\in e}x_i$ such 
that $e\in E(\mathcal{C})$.  
The assignment $\mathcal{C} \mapsto I(\mathcal{C})$ establishes a natural 
one to one correspondence between the family of clutters and the family of 
square-free monomial ideals. A subset $F$ of $X$ is called 
{\it independent\/} or {\it stable\/} if $e\not\subset F$ for any  
$e\in E(\mathcal{C})$. The dual concept of a stable vertex set
is a 
{\it vertex cover\/}, i.e., a subset $C$ of $X$ is a vertex cover of
$\mathcal{C}$  
if and only if $X\setminus C$ is a stable vertex set. 
A first hint of the rich interaction between the combinatorics of $\CC$ 
and the algebra of $I(\CC)$ is that the
number of vertices in a minimum vertex cover of $\mathcal{C}$ 
(the {\it covering number\/} of $\mathcal{C}$) coincides with 
${\rm ht}\, I(\mathcal{C})$, the {\it height\/} of the ideal $I(\mathcal{C})$. 

If $e$ is an edge of $\mathcal{C}$, its {\it
characteristic vector\/} is the vector $v=\sum_{x_i\in e}e_i$, where
$e_i$ is the $i^{\textup{th}}$ unit vector in $\mathbb{R}^n$.
Let $\A = \{ v_1, \ldots, v_q \} \subset \{ 0,1\}^n$ 
denote the characteristic vectors of the edges of $\CC$ and let $A$ denote 
the matrix whose columns, in order, are the vectors of $\A$. We call $A$ 
the {\em clutter matrix\/} or {\it incidence matrix\/} of $\CC$. 
The {\it Ehrhart ring\/} of the lattice polytope 
$P={\rm conv}(\mathcal{A})$ 
is the $K$-subring of $R[t]$  given by
$$
A(P)=K[\{x^at^b\, \vert \, a \in b P \cap\mathbb{Z}^n\}],
$$
where $t$ is a new variable and $bP= \{ b p \, \vert \, p \in P \}$ 
for each $b \in \N$. 
We use  $x^a$ as an abbreviation for $x_1^{a_1} \cdots x_n^{a_n}$, 
where $a=(a_i) \in \mathbb{N}^n$. The {\it homogeneous subring\/}  
of $\mathcal{A}$ is the monomial subring
$$
K[x^{v_1}t,\ldots,x^{v_q}t]\subset R[t].
$$
This ring is in fact a standard graded $K$-algebra because the vector $(v_i,1)$ 
lies in the affine hyperplane with last coordinate equal to $1$ for every $i$. 
In general we have the containment 
\begin{equation} \label{ehrhartContainment}
K[x^{v_1}t,\ldots,x^{v_q}t]\subset A(P),
\end{equation}
but as can be seen in \cite{ehrhart,reesclu},  the algebraic properties of
edge ideals and Ehrhart rings of clutters are more tractable when 
the equality holds in this containment. We call such clutters 
{\em Ehrhart clutters} (or we say that the clutter is Ehrhart).

A finite set $\mathcal{H}\subset\mathbb{Z}^n$ is called a {\it Hilbert basis\/} if
$\mathbb{N}\mathcal{H} = \mathbb{R}_+\mathcal{H} \cap \mathbb{Z}^n$, 
where $\mathbb{R}_+\mathcal{H}$ and $\mathbb{N}{\mathcal H}$ are the 
non-negative real span and non-negative integer span respectively 
of ${\mathcal H}$.
It is not hard to see that $\mathcal{C}$ is an Ehrhart 
clutter if and only if the $q$ vectors $\{(v_1,1),\ldots,(v_q,1)\} \subset \{ 0,1\}^{n+1}$
form a Hilbert basis. 

In this article we present two new families of Ehrhart clutters and we 
then use this information to study some algebraic properties of 
$I(\mathcal{C})$ and $A(P)$, 
such as normality, torsion freeness, Castelnuovo-Mumford
regularity and $a$-invariant. The first two properties for edge ideals have already
have been studied before in \cite{ainv,normali,reesclu,clutters,ITG}.
The Castelnuovo-Mumford regularity (see
Definition~\ref{casteluovo-mumforrd-def}) of a graded algebra is a
numerical invariant that measures the  
``complexity'' of its minimal graded free resolution and plays an
important role in computational commutative algebra \cite{BHer,cca}. 
The $a$-invariant of the Ehrhart ring $A(P)$ is the 
largest integer $a\leq -1$ for which $-aP$ has an interior lattice
point \cite{bruns-gubeladze}. In Section~\ref{cas-mum-reg} we
introduce the  regularity and the $a$-invariant in 
combinatorial and algebraic terms. 

On the other hand, a clutter being Ehrhart will enable us to prove combinatorial 
properties, like when certain clutters have the max-flow min-cut property. 
This property is of central importance in combinatorial optimization 
\cite{cornu-book} and so we define it here: the clutter $\CC$ is said to 
have the {\em max-flow min-cut property} (or we say that $\CC$ is MFMC) if the linear program:
\begin{equation} \label{linearprogram}
\textup{max} \{ \langle{\bf 1},y\rangle \, \vert \, y \geq 0, \, A y \leq w \}
\end{equation} 
has an integral optimal solution for all $w \in \N^n$. 
Here $\langle\ ,\, \rangle$ denotes the standard 
inner product and $\mathbf{1}$ is the vector with all its entries
equal to $1$.

The contents of this paper are as follows. The main theorem in 
Section \ref{cas-mum-reg} is a sharp upper bound for the
Castelnuovo-Mumford regularity of $A(P)$. Before stating the 
theorem, recall that a clutter is called $d$-{\it uniform} if all its
edges have size $d$. A clutter is called {\it unmixed}  if all its
minimal vertex covers have the same size. Unmixed clutters and $d$-uniform 
clutters have been studied in \cite{MRV,unmixed} and \cite{mfmc} respectively.

\medskip

\noindent {\bf Theorem \ref{ehrhart-unmixed}}{\it\ If $\mathcal{C}$ is a
$d$-uniform, unmixed MFMC clutter with covering number $g$, then $\CC$ is Ehrhart, 
the $a$-invariant of $A(P)$ is bounded from above by $-g$, and 
the Castelnuovo-Mumford regularity of $A(P)$ is sharply bounded from
above by $(d-1)(g-1)$.}

\medskip 

A key ingredient to showing this result is a formula of
Danilov-Stanley that expresses the canonical module of $A(P)$ using 
polyhedral geometry (see Eq.~(\ref{march14-09-belleza-32})). For
uniform unmixed MFMC clutters, this formula can be made explicit
enough (see Eq.~(\ref{jan24-09-4})) to allow to prove our estimates
for 
the regularity and
the $a$-invariant of $A(P)$.

The {\it blocker\/} of
a clutter $\mathcal{C}$, denoted by $\Upsilon(\mathcal{C})$, is the clutter 
whose edges are the minimal vertex covers of $\mathcal{C}$ (minimal
with respect to inclusion). Sometimes the {\it blocker\/} of a clutter
is referred to as the {\it Alexander dual\/} of the clutter. The edge
ideal of $\Upsilon(\mathcal{C})$ is called the ideal of {\it vertex
covers} of $\mathcal{C}$ or the {\it Alexander dual\/} of $I(\mathcal{C})$. 
As a corollary of Theorem \ref{ehrhart-unmixed}, using the fact that the 
blocker of a bipartite graph satisfies the max-flow min-cut property \cite{Schr2}, 
we obtain: 

\medskip

\noindent {\bf Corollary \ref{ehrhart-unmixed-coro}}{\it\ Let $G$ be an unmixed
bipartite graph with $n$ vertices, let $\mathcal{A}=\{v_1,\ldots,v_q\}$ be the set
of column vectors of the clutter matrix of the blocker of $G$, and let
$P={\rm conv}(\mathcal{A})$. Then the blocker of $G$ is Ehrhart 
and the Castelnuovo-Mumford regularity of $A(P)$ is bounded by $(n/2)-1$.}

\medskip

In Section \ref{meyniel}, we turn our attention to the clique clutters 
of Meyniel graphs. A {\it clique\/} of a graph is a set of mutually adjacent 
vertices.  The {\it clique clutter\/} of a graph $G$, denoted
by ${\rm cl}(G)$, is the clutter on $V(G)$ whose edges are the 
maximal cliques of $G$. The clutter matrix of ${\rm cl}(G)$ is
called the {\it vertex-clique\/} matrix
of $G$. A {\it Meyniel graph\/} is a simple graph in which
every odd cycle of length at least five has at least two chords, where a {\it
chord} of a cycle $C$ is an edge joining two 
non-adjacent vertices of $C$. A clutter $\mathcal{C}$ is called {\it
ideal\/} if the polyhedron $Q(A)=\{x\vert\, x\geq 0;xA\geq
\mathbf{1}\}$ has only integral vertices, where $A$ is the clutter matrix of
$\mathcal{C}$. Our main result in 
Section~\ref{meyniel} is:

\medskip

\noindent {\bf Theorem \ref{thm:meynielMFMC}}
{\em Let $\CC$ be the clique clutter of a Meyniel graph. If $\CC$ is
ideal, then $\CC$ is MFMC.}

\medskip

Central to proving this result is that the clique clutters of Meyniel graphs 
are Ehrhart, the proof of which arises chiefly from a polyhedral interpretation 
of a known characterization of Meyniel graphs (see Theorem~\ref{meyniel-strongly-perfect}) 
and the fact that the {\em cone} of a vertex over a graph preserves the 
Meyniel property (see Lemma~\ref{lem:suspension}). 
Theorem \ref{thm:meynielMFMC} can also be stated as follows: 
the clique clutter of a Meyniel graph $G$ is ideal if and only if 
$I^i=I^{(i)}$ for $i\geq 1$, where $I\subset R$ is the edge ideal of the 
clique clutter of $G$ and $I^{(i)}$ is the $i^{\textup{th}}$ symbolic
power of $I$. This algebraic 
perspective plays a starring role in the proof of Theorem \ref{thm:meynielMFMC} 
and will be described in great detail in Section~\ref{meyniel}.

Let us take this opportunity to justify the importance of Theorem~\ref{thm:meynielMFMC}. 
Inspired by Lov{\' a}sz's weak perfect graph theorem (see
Theorem~\ref{thm:perfect}), Conforti and Cornu\'ejols 
conjectured \cite[Conjecture 1.6]{cornu-book} that if $\CC$ has the 
{\em packing property} 
(i.e., the linear program (\ref{linearprogram}) has an integer optimal solution 
for all $\omega \in \{0,1, \infty \}^n$), then $\CC$ is also MFMC. 
However, the packing property has proved quite difficult to understand and so, 
given that the Edmonds-Giles theorem \cite[Corollary 22.1c]{Schr} implies that if 
$\CC$ is MFMC then $\CC$ is ideal, 
some energies have been devoted to instead asking: if $\CC$ is an
ideal clutter, then what additional properties on $\CC$ will suffice for $\CC$ to be MFMC? 
For example, one property that suffices is the diadic property \cite[Theorem 1.3]{cgm}. 
We conjecture that the following holds:

\begin{conjecture} \label{con:ideal}
{\it Let $\CC$ be the clique clutter of a perfect graph. If $\CC$ is
ideal, then $\CC$ is MFMC.}
\end{conjecture}

Experimentally, Conjecture~\ref{con:ideal} holds in each of the many distinct examples of 
perfect graphs in \cite[\S 7]{hougardy}, verified using a combination of the 
computational programs {\tt Normaliz} \cite{normaliz2} and {\tt Polymake} \cite{polymake}. 
Since every Meyniel graph is perfect \cite[Theorem~66.6]{Schr2}, then 
Theorem~\ref{thm:meynielMFMC} states that the conjecture holds for 
Meyniel graphs. Conjecture~\ref{con:ideal} also holds when the clique
clutter $\mathcal{C}$ of a perfect graph is uniform \cite[Corollary
2.9]{perfect}. In 
Theorem~\ref{thm:uniformMFMC} we 
provide a simpler alternative proof of the uniform case, again by showing that 
these clutters are Ehrhart.

Section \ref{meyniel} is closed with two examples of clique clutters of perfect graphs. 
The first example shows that the common approach of
Theorem~\ref{thm:meynielMFMC} and Theorem~\ref{thm:uniformMFMC} involving 
Ehrhart clutters is not one that can be relied upon to 
prove Conjecture~\ref{con:ideal} outright. The second example is a perfect graph 
whose clique clutter edge ideal is not normal, in sharp contrast to a central 
result of \cite{perfect} which shows that the
edge ideal of the blocker of a perfect graph is always normal. Thus finding a 
graph theoretical description for the normality of edge ideals of clique clutters 
of perfect graphs remains an open problem.

We close the paper by providing some characterizations of total dual 
integrality, using a generalization of Ehrhart clutters. 
We say that the system $xA\leq w$ is 
{\it totally dual integral\/} (TDI for short) if the minimum in 
the LP-duality equation
$$
{\rm max}\{\langle a,x\rangle \vert\, xA\leq w\}=
{\rm min}\{\langle y,w\rangle \vert\, 
y\geq 0;\, Ay=a \}
$$
has an integral optimum solution $y$ for each integral vector $a$ with 
finite minimum. Note that the MFMC property for a clutter $\CC$ in the 
previous sections can be stated as $x[A|I_n] \leq (-{\bf 1}|0)$ is TDI, 
where $A$ is the clutter matrix of $\CC$, $\mathbf{1}$ is the vector 
of all $1$'s and $I_n$ is an identity matrix.

A rational polyhedron $Q$ is called {\it integral\/} if $Q$ is the 
convex hull of the integral points in $Q$. A classical theorem of 
Edmonds and Giles is that if the system $xA\leq w$ 
is TDI, then the polyhedron $\{x \, \vert \, xA\leq w\}$ is integral 
\cite[Corollary~22,1c]{Schr}. Its converse does not hold in general 
so, similar to Section \ref{meyniel}, it is natural to ask: what properties 
can be added to a matrix $A$ so that $\{x \, \vert \, xA\leq w\}$ being integral 
implies that $xA\leq w$ is TDI? For example, Lovasz's weak perfect graph 
theorem mentioned above can be restated as such a converse holding. We show 
the following theorem:

\medskip

\noindent {\bf Theorem~\ref{viwi-tdi}} {\em Let $A$ be an integral matrix with 
column vectors $v_1,\ldots,v_q$ and let $w=(w_i)$ be an integral vector. 
If the polyhedron $P=\{x\vert\, xA\leq w\}$ is integral 
and $\mathcal{H}(A,w)=\{(v_i,w_i)\}_{i=1}^q$
is a Hilbert 
basis, then the system $xA\leq w$ is TDI.}

\medskip

Note that the set of vectors $\mathcal{H}(A,w)$ being a Hilbert basis is in some 
sense a generalization of Ehrhart clutters. 
We end the section with Proposition~\ref{tdi-non-negative} describing a scenario 
where the converse to Theorem~\ref{viwi-tdi} holds.

\section{Castelnuovo-Mumford regularity and $a$-invariants}\label{cas-mum-reg}

We continue
using the 
definitions and terms from the introduction. 
In this section we give sharp upper bounds for the regularity and the
$a$-invariant of Ehrhart rings arising from uniform unmixed MFMC clutters. 

First we introduce the $a$-invariant and the regularity in
combinatorial and algebraic terms. Assume that 
$A(P)=K[x^{v_1}t,\ldots,x^{v_q}t]$, i.e., assume that $\mathcal{C}$ is
an Ehrhart clutter. Then $A(P)$ becomes a standard graded $K$-algebra
$$
A(P)=\bigoplus_{i=0}^\infty A(P)_i
$$ 
with $i^{\textup{th}}$ component given by
$$
A(P)_i=\sum_{a\in\mathbb{Z}^n\cap iP}K x^a t^i.
$$

A nice property of $A(P)$ is its normality, i.e., $A(P)$ is an
integral domain which is integrally closed 
in its field of fractions \cite[p.~276]{BHer}. Therefore $A(P)$ is a 
Cohen-Macaulay domain by a theorem of Hochster \cite{Ho1}. The {\it Hilbert
series\/} of $A(P)$ is given by 
$$
F(A(P),z)=\sum_{i=0}^\infty\dim_K A(P)_iz^i=\sum_{i=0}^\infty |\mathbb{Z}^n\cap iP|z^i,
$$
this series is called the {\it Ehrhart series\/}
of $P$. By the Hilbert-Serre theorem \cite{BHer,Sta1}, and the fact that
$A(P)$ is a Cohen-Macaulay domain, it follows that this is a 
rational function that can be uniquely written as: 
$$
F(A(P),z)=\frac{h(z)}{(1-z)^{d+1}}=\frac{h_0+h_1z+\cdots+h_sz^s}{(1-z)^{d+1}},
$$
with $h(1)>0$, $h_i\in\mathbb{N}$ for all 
$i$, $h_s>0$ and $d=\dim(P)$. The $a$-{\it invariant\/} of $A(P)$, denoted by
$a(A(P))$, is the degree of $F(A(P),z)$
 as a rational function. This invariant is of combinatorial interest
 because it turns out that $-a(A(P))$ is the 
smallest integer $k\geq 1$ for which $kP$ has an interior lattice
point (see \cite[Theorem~6.51]{bruns-gubeladze}). 

The vector $h=(h_0,\ldots,h_s)$ is called the $h$-{\it vector\/} of $A(P)$.
As $A(P)$ is a Cohen-Macaulay standard  graded $K$-algebra, according 
to \cite[Corollary~B.4.1, p.~347]{Vas1}, the number $s$
turns out to be ${\rm reg}(A(P))$, the Castelnuovo-Mumford regularity of
$A(P)$ (see Definition~\ref{casteluovo-mumforrd-def}). Thus ${\rm
reg}(A(P))$ measures the size of the $h$-vector of $A(P)$ and we have
the equality
$${\rm reg}(A(P))=\dim(A(P))+a(A(P)).$$ 

The $h$-vector of $A(P)$ is of interest in
algebra and combinatorics
\cite{bruns-gubeladze,BHer,hibi-ehrhart,cca,stanley-rational-polytopes}
because it encodes information about 
the lattice polytope $P$ and the algebraic structure of
$A(P)$. For instance $h(1)$ is the multiplicity of the ring 
$A(P)$ and $h(1)=d!{\rm vol}(P)$, where  
${\rm vol}(P)$ is the relative volume of $P$. 

Next we give the definition of regularity of a homogeneous subring 
in terms of its minimal graded free resolution.

\begin{definition}\label{casteluovo-mumforrd-def} 
Let $S=K[x^{v_1}t,\ldots,x^{v_q}t]$ be a homogeneous subring 
with the standard grading induced by
$\deg(x^at^b)=b$. Let
$$
K[t_1,\ldots,t_q]/I_{\mathcal{A}}\simeq S,\ \ \ \ \ \
\overline{t}_i\mapsto x^{v_i}t,
$$
be a presentation of $S$, and let $\mathbb{F}_\star$ be the
minimal graded resolution of $S$ by free $K[t_1,\ldots,t_q]$-modules.
The {\it Castelnuovo-Mumford regularity\/} of $S$ is 
defined as ${\rm reg}(S)=\max\{b_j-j\}$, where $b_j$ is the maximum of
the degrees of a minimal set of generators of $F_j$, the
$j^{\textup{th}}$ component of $\mathbb{F}_\star$.
\end{definition}

\begin{proposition}{\rm \cite[Proposition~5.8]{reesclu}}\label{nov20-03-chordal} 
Let $ \mathcal{C}$ be a
$d$-uniform clutter and let $A$ be its clutter matrix. If the 
polyhedron $Q(A)=\{x\vert\, x\geq 0;xA\geq\mathbf{1}\}$ is integral, then 
there are $X_1,\ldots,X_d$ mutually disjoint minimal vertex covers of
$\mathcal C$ such  that $X=\cup_{i=1}^d X_i$. In particular if
$g_1,\ldots,g_q$ are the edges of $\mathcal{C}$, 
$|X_i \cap g_j| = 1$ for all $i,j$.
\end{proposition}

We come to the main result of this section.

\begin{theorem}\label{ehrhart-unmixed} Let $\mathcal{C}$ be a
$d$-uniform unmixed MFMC clutter with covering number $g$ and let 
$\mathcal{A}=\{v_1,\ldots,v_q\}$ be the characteristic vectors of the
edges of $\mathcal{C}$. If $A(P)$ is the  
Ehrhart ring of $P={\rm conv}(\mathcal{A})$, then $\mathcal{C}$ is an
Ehrhart clutter, the $a$-invariant of $A(P)$ is bounded from above by $-g$ and 
the Castelnuovo-Mumford regularity of $A(P)$ is sharply bounded from
above by
$(d-1)(g-1)$.
\end{theorem}

\begin{proof} Let $\mathcal{B}=\{(v_i,1)\}_{i=1}^q$ and 
$\mathcal{A}'=\mathcal{B}\cup\{e_i\}_{i=1}^n$, where $n$ is the number
of vertices of $\mathcal{C}$ and $e_i$ is the $i^{\textup{th}}$ unit
vector. We first show the equality 
\begin{equation}\label{jan24-09-1}
\mathbb{R}_+\mathcal{B}=\mathbb{R}\,
\mathcal{B}\cap\mathbb{R}_+\mathcal{A}',
\end{equation}
where $\mathbb{R}\,\mathcal{B}$ is the vector space spanned by
$\mathcal{B}$ and $\mathbb{R}_+\mathcal{B}$ is the cone generated by
$\mathcal{B}$. The left hand side is clearly contained in the right hand side.
Conversely, take $(a,b)$ in the cone $\mathbb{R}\,
\mathcal{B}\cap\mathbb{R}_+\mathcal{A}'$, where $a\in\mathbb{R}^n$ 
and $b\in\mathbb{R}$. Then one has
\begin{eqnarray*}
(a,b)&=&\eta_1(v_1,1)+\cdots+\eta_q(v_q,1) \ \ \ \ \ (\eta_i\in
\mathbb{R}),\\
(a,b)     &=&
     \lambda_1(v_1,1)+\cdots+\lambda_q(v_q,1)+\mu_1e_1+\cdots+\mu_ne_n\
     \ \ \ \ (\lambda_i,\, \mu_j\in\mathbb{R}_+\, \forall\, i,j). 
\end{eqnarray*}
For $a=(a_i)\in\mathbb{R}^n$, we set $|a|=\sum_{i}a_i$. 
Hence using that $\mathcal{C}$ is $d$-uniform, i.e., $|v_i|=d$ for all $i$, we get 
$bd=bd+\sum_{i}\mu_i$. This proves that $\mu_i=0$ for all $i$ 
and thus $(a,b)$ is in $\mathbb{R}_+\mathcal{B}$, as required.

Next we prove that $\mathcal{C}$ is an Ehrhart clutter, i.e., we will prove 
the equality 
\begin{equation}\label{jan24-09}
K[x^{v_1}t,\ldots,x^{v_q}t]=A(P).
\end{equation}
By \cite[Theorem~4.6]{reesclu}, 
the Rees algebra 
$$
R[I(\mathcal{C})t]=R[x^{v_1}t,\ldots,x^{v_q}t]\subset R[t]
$$
of the edge ideal $I(\mathcal{C})=(x^{v_1},\ldots,x^{v_q})$ is
normal. Hence, using \cite[Theorem~3.15]{ehrhart}, we obtain the required equality.

The next step in the proof is to find a good expression for 
the canonical module of $A(P)$ (see Eq.~(\ref{jan24-09-4}) below) 
that can be used to estimate the regularity 
and the $a$-invariant of $A(P)$. We begin by extracting some of the
information encoded in the polyhedral representation of the 
cone $\mathbb{R}_+\mathcal{A}'$. Let $C_1,\ldots,C_s$ be the minimal vertex covers
of $\mathcal{C}$ and let $u_k=\sum_{x_i\in C_k}e_i$ for $1\leq k\leq
s$. By
\cite[Proposition~3.13 and Theorem~4.6]{reesclu} we obtain that the 
irreducible representation of $\mathbb{R}_+\mathcal{A}'$ as an
intersection of closed halfspaces is given by  
\begin{equation}\label{jan24-09-2}
\mathbb{R}_+\mathcal{A}'=H_{e_1}^+\cap\cdots\cap H_{e_{n+1}}^+\cap
H_{(u_1,-1)}^+\cap\cdots\cap H_{(u_s,-1)}^+.
\end{equation}
Here $H_{a}^+$ denotes the closed halfspace 
$H_a^+=\{x\vert\, \langle
x,a\rangle\geq 0\}
$
and $H_a$ stands for the hyperplane through the origin with normal
vector $a$. Let $A$ be the clutter matrix of $\mathcal{C}$ whose
columns are $v_1,\ldots,v_q$. The set covering polyhedron 
$$
Q(A)=\{x\vert\, x\geq 0; xA\geq\mathbf{1}\}
$$
is
integral \cite[Theorem~4.6]{reesclu} and $\mathcal{C}$ is unmixed by
hypothesis. Therefore, by Proposition~\ref{nov20-03-chordal}, there
are $X_1,\ldots,X_d$ mutually disjoint minimal vertex covers of
$\mathcal C$ of size $g$ such  that $X=\cup_{i=1}^d X_i$. 
Notice that 
$|X_i\cap f|=1$ for $1\leq i\leq d$ and $f\in E(\mathcal{C})$. We may
assume that $X_i=C_i$ for $1\leq i\leq d$. 
Therefore, using Eqs.~(\ref{jan24-09-1}) and
(\ref{jan24-09-2}), we get 
\begin{eqnarray}
\mathbb{R}_+\mathcal{B}&=&\mathbb{R}\,
\mathcal{B}\cap\mathbb{R}_+\mathcal{A}'\nonumber\\
&=&\mathbb{R}\, \mathcal{B}\cap H_{e_1}^+\cap\cdots\cap H_{e_{n+1}}^+\cap
H_{(u_1,-1)}^+\cap\cdots\cap H_{(u_s,-1)}^+\nonumber\\
&=&\mathbb{R}\, \mathcal{B}\cap H_{e_1}^+\cap\cdots\cap
H_{e_{n}}^+\cap H_{e_{n+1}}^+\cap\left(\cap_{i\in
\mathcal{I}}H_{(u_i,-1)}^+\right),\label{jan24-09-3}
\end{eqnarray}
where $i\in \mathcal{I}$ if and only if $H_{(u_i,-1)}^+$ defines a proper 
face of the cone $\mathbb{R}_+\mathcal{B}$. 
As $(v_i,1)$ lies in the affine hyperplane 
$x_{n+1}=1$ for all $i$, the ring $A(P)$ becomes a graded $K$-algebra generated
by monomials of degree $1$. 
Notice that a monomial $x^at^b$ has degree $b$ in this grading. The
Ehrhart ring $A(P)$ is a normal domain. Then, according to a well 
known formula of  Danilov-Stanley \cite[Theorem 6.3.5]{BHer},  its canonical module
is the ideal of $A(P)$ given by   
\begin{equation}\label{march14-09-belleza-32}
\omega_{A(P)}=(\{x_1^{a_1}\cdots x_n^{a_n}   
t^{a_{n+1}}\vert\, a=(a_i)   
\in\mathbb{N}\mathcal{B}\cap(\mathbb{R}_+\mathcal{B})^{\rm o}\}),   
\end{equation}
where $(\mathbb{R}_+\mathcal{B})^{\rm o}$ denotes the relative
interior of the cone $\mathbb{R}_+\mathcal{B}$. Using
Eqs.~(\ref{jan24-09}) and (\ref{jan24-09-3}) we can express the
canonical module as:
\begin{equation}\label{jan24-09-4}   
\omega_{A(P)}=(\{x_1^{a_1}\cdots x_n^{a_n}   
t^{a_{n+1}}\vert\, a=(a_i)   
\in\mathbb{R}\, \mathcal{B};\, a_i\geq 1\,\forall\, i; \textstyle\sum_{x_i\in
C_k}a_i\geq a_{n+1}+1\ \mbox{ for }\ k\in\mathcal{I}\}). 
\end{equation}   
Next we estimate the $a$-invariant of $A(P)$. 
Recall that the $a$-{\it invariant\/} of
$A(P)$ is  the degree, as a rational 
function, of the Hilbert series of $A(P)$ \cite[p.~99]{monalg}. The
ring $A(P)$ is 
normal, then $A(P)$ is Cohen-Macaulay \cite{Ho1} and its $a$-invariant is
given by
\begin{equation}\label{princeton-fall-07-1}
a(A(P))=-{\rm min}\{\, i\, \vert\, (\omega_{A(P)})_i\neq 0\},
\end{equation}
see \cite[p.~141]{BHer} and \cite[Proposition~4.2.3]{monalg}. Take an    
arbitrary monomial    
$x^at^b=x_1^{a_1}\cdots x_n^{a_n}t^b$ in the ideal $\omega_{A(P)}$. By
Eqs.~(\ref{jan24-09-3}) and (\ref{jan24-09-4}), the vector $(a,b)$ is
in $\mathbb{R}_+\mathcal{B}$ and $a_i\geq 1$ 
for all $i$. Thus we can write
$$
(a,b)=\lambda_1(v_1,1)+\cdots+\lambda_q(v_q,1)\ \ \ \ (\lambda_i\geq
0).
$$
Since $\langle v_i,u_k \rangle=1$ for $i=1,\ldots,q$ and
$k=1,\ldots,d$, we obtain
$$
g=|u_k|\leq\sum_{x_i\in C_k}a_i=\langle
a,u_k\rangle=\lambda_1\langle v_1,u_k\rangle+\cdots+
\lambda_q\langle v_q,u_k\rangle=\lambda_1+\cdots+\lambda_q=b
$$
for $1\leq k\leq d$. This means that $\deg(x^at^b)\geq g$.
Consequently $-a(A(P))\geq g$,
as required. Next we show that ${\rm reg}(A(P))\leq(d-1)(g-1)$. Since $A(P)$ is Cohen-Macaulay, we
have 
\begin{equation}\label{jan25-09}
{\rm reg}(A(P))=\dim(A(P))+a(A(P))\leq \dim(A(P))-g,
\end{equation}
see \cite[Corollary~B.4.1, p.~347]{Vas1}. Using that 
$\langle v_i,u_k \rangle=1$ for $i=1,\ldots,q$ and
$k=1,\ldots,d$, by induction on $d$ it is seen that 
${\rm rank}(A)\leq g+(d-1)(g-1)$. Thus using the fact that
$\dim(A(P))={\rm rank}(A)$ and Eq.~(\ref{jan25-09}), we 
get ${\rm reg}(A(P))\leq (d-1)(g-1)$. 

Finally, we now show that the upper bounds for the $a$-invariant and
for the regularity are sharp. Let
$\mathcal{C}$ be the clutter with vertex set $X=\cup_{i=1}^dX_i$ 
whose minimal vertex covers are exactly $X_1,\ldots,X_d$. Let
$v_1,\ldots,v_q$ be the characteristic vectors of the edges of
$\mathcal{C}$ and let $A$ be the matrix with column vectors
$v_1,\ldots,v_q$. Using \cite[Corollary~83.1a]{Schr2} (cf.
\cite[Corollary~4.26]{reesclu}) it is not hard
to see that $\mathcal{C}$ satisfies the hypotheses of the theorem,
i.e., the clutter $\mathcal{C}$ is MFMC, 
is $d$-uniform, unmixed and has covering number equal to $g$. 
Moreover the rank of $A$ is $g+(d-1)(g-1)$. Thus by
Eq.~(\ref{jan25-09}) it suffices to show that $a(A(P))=-g$. 
Any edge of $\mathcal{C}$ intersects any minimal vertex
cover of $\mathcal{C}$ in
exactly one vertex. Therefore, using Eq.~(\ref{jan24-09-3}), we get 
\begin{equation}
\mathbb{R}_+\mathcal{B}=\mathbb{R}\, \mathcal{B}\cap H_{e_1}^+\cap\cdots\cap
H_{e_{n}}^+\cap H_{e_{n+1}}^+.
\end{equation}
Hence, using Eq.~(\ref{jan24-09-4}), we can express the
canonical module as:
\begin{equation}\label{jan29-09}   
\omega_{A(P)}=(\{x_1^{a_1}\cdots x_n^{a_n}   
t^{a_{n+1}}\vert\, a=(a_i)   
\in\mathbb{R}\, \mathcal{B};\, a_i\geq 1\mbox{ for }i=1,\ldots,n+1\}). 
\end{equation}   
It is well known that MFMC clutters have the K\"onig property (the
covering number equals the maximum number of mutually disjoint edges). Thus
$\mathcal{C}$ has $g$ mutually disjoint edges whose union is $X$, 
by relabeling the $v_i$'s if
necessary, we may assume that $v_1,\ldots,v_g$ satisfy $\mathbf{1}=v_1+\cdots+v_g$. Thus by
Eq.~(\ref{jan29-09}), we get that the monomial $x_1\cdots
x_nt^g$ belongs to $\omega_{A(P)}$. Consequently $a(A(P))\geq -g$ and
the equality $a(A(P))=-g$ follows. 
\end{proof}

\begin{corollary}\label{ehrhart-unmixed-coro} Let $G$ be an unmixed
bipartite graph with $n$ vertices, let $\mathcal{A}=\{v_1,\ldots,v_q\}$ be the set
of column vectors of the clutter matrix of the blocker of $G$, and
let $P={\rm conv}(\mathcal{A})$. Then the blocker of $G$ is an Ehrhart
clutter and the Castelnuovo-Mumford regularity of $A(P)$ is sharply
bounded from above by $(n/2)-1$. 
\end{corollary}

\begin{proof} Let $\mathcal{C}=\Upsilon(G)$ be the clutter of minimal vertex
covers of the bipartite graph $G$ and let $A$ be the matrix with
column vectors $v_1,\ldots,v_q$. Since $A$ is the clutter matrix of
$\mathcal{C}$ and all cycles of $G$ are even, 
it is well known \cite[Theorem 83.1a(v)]{Schr2} 
that the clutter $\mathcal{C}$ has the max-flow min-cut property. 
The covering number of
$\mathcal{C}$ is equal to $2$ because the blocker of $\mathcal{C}$ is
$G$. Moreover, as $G$ is bipartite and has no isolated
vertices, it is seen that $n$ is even and that all edges of
$\mathcal{C}$ have size $n/2$ (see for instance
\cite[Lemma~6.4.2]{monalg}). Therefore by
Theorem~\ref{ehrhart-unmixed}, the Castelnuovo-Mumford regularity of $A(P)$ is bounded by
$(n/2)-1$. 
\end{proof}

\section{Clique clutters with the Ehrhart Property} \label{meyniel}

The main result of this section is that Conjecture~\ref{con:ideal} holds for Meyniel graphs.
\begin{theorem} \label{thm:meynielMFMC}
Let $\CC$ be the clique clutter of a Meyniel graph. If $\CC$ is ideal, then $\CC$ is MFMC.
\end{theorem}
We prove this result by studying the algebraic properties of edge ideals of clutters 
and by showing that clique clutters of Meyniel graphs are Ehrhart. As noted in the 
introduction, Conjecture~\ref{con:ideal} also holds for clique clutters of perfect graphs 
that are $d$-uniform (all edges have cardinality equal to $d$). We present a new simpler 
proof of that statement here, the heart of which is the same as the proof in the 
case of Meyniel graphs. Finally, we finish with examples of perfect graphs whose 
clique clutters are not Ehrhart, thus showing that a different approach than that 
presented here is needed to completely resolve Conjecture~\ref{con:ideal}. 

We begin with the necessary algebraic background. Let $\mathcal{C}$ be
any clutter and let $C_1,\ldots,C_s$ be the minimal vertex covers of
$\mathcal{C}$. By \cite[Proposition 6.1.16]{monalg}, 
the primary decomposition of the edge ideal of $\CC$ is given by 
$$
I(\mathcal{C})=\mathfrak{p}_1\cap
\cdots\cap \mathfrak{p}_s,
$$
where $\mathfrak{p}_i= (C_i)$ for $1\leq i\leq s$ and $(C_i)$ denotes the 
prime ideal of $R$ generated by the minimal vertex cover $C_i$. 
The $i^{\textup{th}}$ {\it symbolic power\/} of
$I=I(\mathcal{C})$ is the ideal of $R$ given by
$$
I^{(i)}= \mathfrak{p}_1^i\cap\cdots\cap \mathfrak{p}_s^i,
$$
and the integral closure of $I^i$ is the ideal of $R$  
given by (see \cite{monalg}):
$$
\overline{I^i}=(\{x^a\in R\vert\, \exists\, p\geq
1\ \mbox{such that}\ (x^a)^{p}\in I^{pi}\}).
$$ 
A central result in this area shows that a clutter $\CC$ is MFMC 
if and only if its edge ideal $I$ is {\em normally torsion free}, 
i.e., if and only if $I^i=I^{(i)}$ for $i\geq 1$ \cite{clutters}. 
The proof of the following result is essentially the same as that made 
in \cite[Corollary 2.9]{perfect}.
\begin{theorem} \label{theorem-ehrhart-mfmc} 
Let $\CC$ be a clutter. If $\CC$ is both Ehrhart and ideal, then $\CC$ is MFMC.
\end{theorem}
\begin{proof} Let $\{v_1,\ldots,v_q\}$ be the set of columns of the clutter
matrix of $\mathcal{C}$ and let $I=I(\CC)$ be the edge ideal of 
$\mathcal{C}$. Assuming that 
$\CC$ is an Ehrhart clutter, we show that the following four conditions 
are equivalent{\rm:} 
\begin{enumerate}
\item[\rm (i)] $\CC$ is MFMC.
\item[\rm (ii)] $I^i=I^{(i)}$ for $i\geq 1$.
\item[\rm (iii)] $\overline{I^i}=I^{(i)}$ for $i\geq 1$.
\item[\rm (iv)] $\CC$ is ideal.
\end{enumerate}
(i) $\Rightarrow$ (ii): This was shown in \cite[Corollary~3.14]{clutters}. 
(ii) $\Rightarrow$ (iii): Follows readily because in general one has 
the inclusions $I^i\subset \overline{I^i}\subset I^{(i)}$ for $i\geq 1$. 
(iii) $\Rightarrow$ (iv): This was shown in \cite[Corollary~3.13]{clutters}. 
(iv) $\Rightarrow$ (i): By our hypothesis $\CC$ is Ehrhart, 
which is equivalent to saying that 
$
K[x^{v_1}t,\ldots,x^{v_q}t]=A(P),
$
where $P={\rm conv}(v_1,\ldots,v_q)$. As $\CC$ is ideal, a direct
application of \cite[Proposition~4.4 and Theorem~4.6]{reesclu} 
gives that the clutter $\CC$ is MFMC. Thus the four conditions are
equivalent. Hence if $\mathcal{C}$ is Ehrhart and ideal, then $\mathcal{C}$ is
MFMC.
\end{proof}

To prove Theorem~\ref{thm:meynielMFMC} we will show that the clique clutters of 
Meyniel graphs are Ehrhart and then apply Theorem~\ref{theorem-ehrhart-mfmc}. 
A starting point to proving Conjecture~\ref{con:ideal} for Meyniel graphs is 
the following classification:

\begin{theorem}{\rm (Ho{\' a}ng \cite{hoang})} \label{meyniel-strongly-perfect} 
A graph $G$ is Meyniel 
if and only if for each induced subgraph $H$ and for each vertex $u$
of $H$, there exists a stable set in $H$ that contains $u$ and 
this set intersects all maximal cliques of $H$.
\end{theorem}
Using this classification we will show that if $\CC$ is the clique clutter of a 
Meyniel graph such that $\{v_1,\ldots,v_q\}$ are the characteristic
vectors of the edges of $\mathcal{C}$, then $\{ (v_1,1),\ldots,(v_q,1)\}$ 
is a Hilbert basis. By Theorem~\ref{theorem-ehrhart-mfmc} we will then have a proof 
of Theorem~\ref{thm:meynielMFMC}. 

We begin by recalling the polyhedral weak perfect graph theorem. Here the {\em complement} 
of a graph $G$ on the vertex set $\{ x_1, \ldots, x_n \}$ is the
graph $G^c$ on the same 
vertices as $G$ but whose maximal cliques are precisely the maximal
stable sets of $G$. 
\begin{theorem} \label{thm:perfect} 
{\rm (Lov{\' a}sz \cite{lovasz}, Chv{\' a}tal \cite{chvatal}, Fulkerson \cite{fulkerson})} 
A graph $G$ is perfect 
if and only if 
its complement is perfect 
if and only if 
the stability polytope 
$\textup{Stab}(G):= \{x\vert\,x \geq 0;xA \leq {\bf 1} \}$
is integral
if and only if 
the system $x \geq 0 \, ; \, xA \leq {\bf 1}$ is TDI, where 
$A$ is the clique clutter matrix of the graph $G$.
\end{theorem}

In this case the vertices of the stability polytope of $G$ are 
precisely the characteristic vectors of the stable sets of $G$. 
Finally, a theorem tying TDI with Hilbert bases: without loss of 
generality, every system of inequalities can be rewritten in the 
form $xA\leq w$. If $F$ is a face of the polyhedron 
$\{ x\vert\, xA\leq w\}$ we say that a column of $A$ is 
{\em active in $F$} if the corresponding inequality in 
$xA\leq w$ is satisfied with equality for all vectors in $F$.
\begin{theorem} \label{thm:tdiHilbert} \cite[Theorem 22.5]{Schr}
The system $xA\leq w$ is TDI if and only if for each face $F$ of the 
polyhedron $P=\{ x\vert\, xA\leq w\}$, the columns of $A$ which are active in $F$ form 
a Hilbert basis.
\end{theorem}

By \cite[Theorem~66.6]{Schr2}, every Meyniel graph is perfect and 
so we can put Theorem~\ref{thm:perfect} and Theorem~\ref{thm:tdiHilbert} 
to good use for Meyniel graphs.

\begin{lemma} \label{lem:meynielHilbert}
Let $G$ be a
Meyniel graph and let $\mathcal{A}=\{v_1,\ldots,v_q\}$ 
be the set of columns of the vertex-clique matrix of $G$. Then $\A =
\{v_1, \ldots, v_q\}$ is a Hilbert basis. 
\end{lemma}
\begin{proof}
By Theorem~\ref{meyniel-strongly-perfect}, for each $x_k$ in 
$V(G)=\{x_1,\ldots,x_n\}$ there 
exists a stable set $B_k$ of $G$ containing $x_k$ and
intersecting all maximal cliques of $G$. Let $\beta_k=\sum_{x_i\in
B_k}e_i$ be the characteristic vector of $B_k$ for $1\leq k\leq n$. 
Note that in general a clique of $G$ and a stable set of $G$ can meet 
in at most one vertex. Then for each $k=1,2,\ldots, n$ we have 
$\langle e_k,\beta_k\rangle=1$ and 
$\langle v_j,\beta_k\rangle=1$ for all $j$. 
Next, let $\beta := \frac{1}{n}\sum_{k=1}^n{\beta_k}$. Note that 
$\beta$ also has the property that 
\begin{equation}\label{march12-09}
\langle v_j,\beta \rangle =1\ \mbox{ and }\
\langle e_i,\beta \rangle>0,\ \ \ \ \ \forall \, i,j.
\end{equation}
Hence, $\beta$ is in the common intersection of all faces of 
${\rm stab}(G)$ of the form $\langle v_j, x \rangle =1$. Call this
intersection $F_\beta$. Thus $\beta$ belongs to the face $F_\beta$ of  
the stability polytope of $G$, rewritten as 
$\{x:x[A |-I_n] \leq ({\bf 1}| 0)\}$, where $I_n$ is the $n \times n$ 
identity matrix and $({\bf 1}| 0)$ is the vector with $q$ $1$'s followed 
by $n$ $0$'s. By \cite[Theorem~66.6]{Schr2}, every Meyniel graph is perfect and 
so by Theorem~\ref{thm:perfect} the system $x[A| -I_n] \leq ({\bf 1}, 0)$ 
is TDI. Now by Eq.~(\ref{march12-09}), the columns of $[A|-I_n]$ which are active in 
$F_\beta$ are precisely the columns of $A$. By Theorem
\ref{thm:tdiHilbert}, the columns of $[A|-I_n]$ which are active in
$F_\beta$ form a Hilbert basis, 
i.e., $\{v_1, \ldots, v_q \}$ is a Hilbert basis as claimed.
\end{proof}

Next, the {\em suspension of a vertex over $\CC$} is the clutter $\CC^+$ with 
vertex set $X\cup\{x_{n+1}\}$ and edge set 
$\{ e \cup \{ x_{n+1} \} \,\vert\,  e \in E(\CC) \}$, where 
$X=\{x_1,\ldots,x_n\}$ is the vertex set of $\mathcal{C}$. The {\em cone\/}
$C(G)$,  over a graph $G$, is the graph obtained by adding a new
vertex $x_{n+1}$ to $G$ and joining every vertex of $G$ to $x_{n+1}$.
If $G$ is a graph and ${\rm cl}(G)$ is its clique clutter, these two
construction 
are related by ${\rm cl}(C(G))={\rm cl}(G)^+$.
The following is a simple lemma showing that the cone operation 
preserves both Meyniel-ness and perfection in graphs, and the
suspension operation  preserves idealness in a clutter.
\begin{lemma}\label{lem:suspension} Let $G$ be a graph and let
$\mathcal{C}$ be a clutter on the vertex set $X$. Then: 
{\rm (i)} $G$ is Meyniel if and only if $C(G)$ is Meyniel. 
{\rm (ii)} $G$ is perfect if and only if $C(G)$ is perfect. 
{\rm (iii)} $\CC$ is ideal if and only if $\CC^+$ is ideal. 
\end{lemma}
\begin{proof}
(i) This follows immediately from the classification for Meyniel graphs in 
Theorem~\ref{meyniel-strongly-perfect} and the construction of $C(G)$. 

(ii) This follows from the complement of the cone over $G$, $C(G)^c$ equals 
the complement of $G$ with the isolated vertex $\{x_{n+1}\}$ and from a graph 
being perfect if and only if its complement is perfect 
(see Theorem~\ref{thm:perfect}).

(iii) This requires a little more work but the proof is still straightforward. 
Lehman \cite[Theorem 1.17]{cornu-book} showed that a clutter is ideal if and only if 
its blocker (defined in the introduction) is ideal. We now show 
that the clutter $\CC$ is ideal if and only if $\CC^+$ is ideal. 

Consider the blocker $\Upsilon(\CC^+)$ of the clutter $\CC^+$. The
edges of $\Upsilon(\CC^+)$ are precisely the  
edges of $\Upsilon(\CC)$ with the solitary edge $\{ x_{n+1} \}$ (viewing the vertex as an edge). 
Let $B$ (resp. $B^+$) denote the clutter matrix for the blocker of
$\CC$ (resp. $\CC^+$). 
By the decomposition theorem for polyhedra \cite[\S 8.9]{Schr}, 
the polyhedron $Q(B^+)=\{x\vert\, x\geq 0;xB^+\geq\mathbf{1}\}$ will simply be 
$$
(\textup{conv}((w_1,0), \ldots, (w_t,0)) + \R^{n+1}_+) \cap 
\{ y \in \R^{n+1}_+: y_{n+1} \geq 1 \},
$$
where $w_1, \ldots, w_t \in \R^n$ are the vertices of $Q(B)$. 
But the intersection above is simply 
$$
\textup{conv}((w_1,1), \ldots, (w_t,1)) + \R^{n+1}_+.
$$
Hence with this formulation we see that $\CC$ is ideal if and only if $\Upsilon(\CC)$ is ideal 
if and only if $w_1, \ldots, w_t \in \R^n$ are integral if and only if $(w_1,1), \ldots, (w_t,1)$ 
are integral if and only if $\Upsilon(\CC^+)$ is ideal if and only if $\CC^+$ is ideal.
\end{proof}

Note that if $\{ v_1, \ldots, v_q \} \subset \{ 0,1 \}^n$ are the characteristic 
vectors of the edges of the clutter $\CC$, then 
$
\{ (v_1,1), \ldots, (v_q,1) \} \subset \{ 0,1 \}^{n+1}
$ 
are the characteristic vectors of the edges of the clutter $\CC^+$. 
With Lemma \ref{lem:suspension}(i) we are now ready prove our main result:

\medskip

\noindent {\em Proof of Theorem~\ref{thm:meynielMFMC}:} 
If $G$ is a Meyniel graph, then so is $C(G)$, the cone over the graph
$G$. The clutter 
matrix of ${\rm cl}(G)^+$ is the matrix $A^+$ whose columns are precisely 
$\{(v_1,1), \ldots, (v_q,1)\}$. Now applying 
Lemma~\ref{lem:meynielHilbert} to the Meyniel graph 
$C(G)$ and noticing the equality ${\rm cl}(C(G))={\rm cl}(G)^+$, we
get that $\{(v_1,1), \ldots, (v_q,1)\}$ is a Hilbert basis, i.e.,
${\rm cl}(G)$ is an Ehrahrt clutter. In addition, if we assume that the clique 
clutter of $G$ is ideal, then by Theorem~\ref{theorem-ehrhart-mfmc} 
the clique clutter must also be MFMC. \hfill $\square$

\medskip

Using arguments very similar to the above for Meyniel graphs, we now give 
a simpler proof of Conjecture \ref{con:ideal} for uniform clutters which 
was first proved in \cite[Corollary 2.9]{perfect}.

\begin{theorem} \label{thm:uniformMFMC}
Let $\CC$ be the clique clutter of a perfect graph with edges $\{v_1, \ldots, v_q \}$ 
that is both ideal and uniform. Then 
{\rm (a)} $\{v_1, \ldots, v_q \}$ is a Hilbert basis. 
{\rm (b)} $\CC$ is Ehrhart. 
{\rm (c)} $\CC$ is MFMC.
\end{theorem}
\begin{proof}
Recall that a clutter $\CC$ is $d$-uniform if all edges have the same 
cardinality $d$. Let $g_1,\ldots,g_q$ be the edges of $\mathcal{C}$. 
From Proposition~\ref{nov20-03-chordal} we have 
that if $\CC$ is ideal, then there are mutually disjoint sets $X_1, \ldots,
X_d$ whose union is $X$ and such that $|g_j \cap X_i| =1$ for each $j
= 1,\ldots, q$ and 
each $i=1,\ldots,d$. A simple proof of this can be seen in 
\cite[Proposition 2.2]{mfmc}. If $\CC$ is also a clique clutter of a perfect 
graph $G$, then each of the $X_i$'s are also maximal stable sets for the graph $G$. 

(a) Let $\chi_i$ represent the $0/1$ characteristic vector of $X_i$ for every 
$i = 1,\ldots, d$ and 
$$
\gamma := \frac{1}{d}(\chi_1 + \cdots + \chi_d) = (\frac{1}{d},
\ldots, \frac{1}{d}).
$$
Just as in the case of the constructed $\beta$ for Meyniel graphs, 
$\gamma$ also has the similar property that 
$$
\langle v_j,\gamma \rangle =1\ \mbox{ and }\
\langle e_i,\gamma \rangle =1/d>0,\ \ \ \ \ \forall \, i,j.
$$ 
and, by the exact same argument as that for $\beta$ in Lemma~\ref{lem:meynielHilbert}, 
$\gamma$ belongs to a face $F_\gamma$ of 
the stability polytope of $G$, and the columns of $[A|-I_n]$ that are
active in $F_\gamma$ are precisely 
the columns of $A$ and they form a Hilbert basis.

(b) By Lemma \ref{lem:suspension} (parts (ii) and (iii)), cones
preserve 
perfection and suspensions preserve idealness res\-pectively. Since
$\CC={\rm cl}(G)$ is $d$-uniform, then 
$\CC^+={\rm cl}(C(G))$ is $(d+1)$-uniform. Hence part (a) can be
applied to the 
clutter $\CC^+$ 
and so $\{(v_1,1), \ldots, (v_q,1) \}$ forms a Hilbert basis.

(c) Given part (b), its a simple consequence of Theorem
\ref{theorem-ehrhart-mfmc}. 
\end{proof}

\begin{remark}
The clutter $\CC$ being Ehrhart is central to both arguments and so classifying other 
clique clutters that are Ehrhart is certainly of interest. The clique clutter of the 
perfect graph in Example~\ref{pepe-contraejemplo} below is not Ehrhart yet it is ideal 
and MFMC (confirmed computationally). Hence, to completely resolve Conjecture~\ref{con:ideal} 
an approach that differs from the one presented here is needed.
\end{remark}

\begin{example}\label{pepe-contraejemplo} 
Let $\mathcal{K}_{2,4}$ be the complete bipartite
graph with vertex set $X=V_1\cup V_2$ and 
bipartition $V_1=\{x_1,x_2\}$, $V_2=\{x_3,x_4,x_5,x_6\}$. The matrix 
$$
\left(\begin{matrix}
1 &1 &1 &1 &0 &0 &0 &0\\
0 &0 &0 &0 &1 &1 &1 &1\\
1 &0 &0 &0 &1 &0 &0 &0\\
0 &1 &0 &0 &0 &1 &0 &0\\
0 &0 &1 &0 &0 &0 &1 &0\\
0 &0 &0 &1 &0 &0 &0 &1
\end{matrix}\right)
$$
is the incidence matrix of $\mathcal{K}_{2,4}$. Thus
this matrix is totally unimodular and its rows are the maximal cliques
of some perfect graph $G$ \cite[Theorem~82.4]{Schr2}, actually $G$ is
the line graph of the bipartite graph $\mathcal{K}_{2,4}$: 

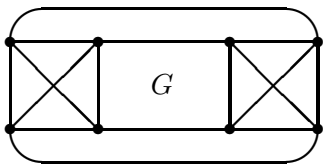
\setlength{\unitlength}{0.023in}
\begin{center}
\begin{picture}(50,20)(12,10)
 \thicklines
 \put(0,0){\line(1,0){70}}
 \put(0,20){\line(1,0){70}}
 \put(0,0){\line(0,0){20}}
 \put(20,0){\line(0,0){20}}
 \put(50,0){\line(0,0){20}}
 \put(70,0){\line(0,0){20}}
 \put(0,0){\line(1,1){20}}
 \put(0,20){\line(1,-1){20}}
 \put(50,0){\line(1,1){20}}
 \put(50,20){\line(1,-1){20}}

 \put(32,8){$G$}

 \put(35,20){\oval(70,15)[t]}
 \put(35,0){\oval(70,15)[b]}

 \put(0,0){\circle*{2.2}}
 \put(0,20){\circle*{2.2}}
 \put(20,0){\circle*{2.2}}
 \put(20,20){\circle*{2.2}}
 \put(50,0){\circle*{2.2}}
 \put(50,20){\circle*{2.2}}
 \put(70,0){\circle*{2.2}}
 \put(70,20){\circle*{2.2}}
\end{picture}
\end{center}
\vspace{1.5cm}

Using {\it Normaliz\/} \cite{normaliz2} it
is seen that $K[x^{v_1}t,\ldots,x^{v_6}t]\subsetneq A(P)$, where
$v_1,\ldots,v_6$ are the rows of the matrix $B$. The extra element 
in $A(P)$ but not in $K[x^{v_1}t,\ldots,x^{v_6}t]$ is the monomial 
with exponent vector $(1,1,1,1,1,1,3)$.
\end{example}

The next example was constructed so that the clique clutter of 
the graph $G$ has a minor (in the sense of hypergraph theory 
\cite{cornu-book}) whose edge ideal is not normal. 

\begin{example}\label{vila-ce} Consider the following graph $G$ with $13$ vertices:
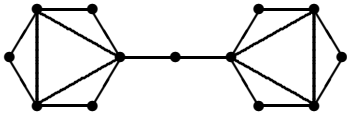
\setlength{\unitlength}{0.029in}
\begin{center}
\begin{picture}(5,5)(19,10)
 \thicklines
 \put(-10,0){\line(3,5){5}}
 \put(10,0){\line(1,0){10}}
 \put(10,0){\line(-3,5){5}}
 \put(-5,-8.66){\line(-3,5){5}}
 \put(5,-8.66){\line(3,5){5}}
 \put(-5,8.66){\line(1,0){10}}
 \put(-5,-8.66){\line(1,0){10}}
 \put(20,0){\line(1,0){10}}

 \put(30,0){\line(3,5){5}}
 \put(50.2,0){\line(-3,5){5}}
 \put(35,-8.66){\line(-3,5){5}}
 \put(35,-8.66){\line(1,0){10}}
 \put(45,-8.66){\line(3,5){5}}
 \put(35,8.66){\line(1,0){10}}

 \put(-10,0){\circle*{2}}
 \put(10,0){\circle*{2}}
 \put(-5,-8.66){\circle*{2}}
 \put(5,8.66){\circle*{2}}
 \put(5,-8.66){\circle*{2}}
 \put(-5,8.66){\circle*{2}}
 \put(20,0){\circle*{2}}

 \put(30,0){\circle*{2}}
 \put(50,0){\circle*{2}}
 \put(35,-8.66){\circle*{2}}
 \put(45,8.66){\circle*{2}}
 \put(45,-8.66){\circle*{2}}
 \put(35,8.66){\circle*{2}}
 \multiput(45,-8.66)(0,.112){150}{\circle*{.02}}
 \multiput(-5,-8.66)(0,.112){150}{\circle*{.01}}
 \multiput(10,0)(-.2,.112){75}{\circle*{.01}}
 \multiput(10,0)(-.2,-.112){75}{\circle*{.01}}
 \multiput(30,0)(.2,.112){75}{\circle*{.02}}
 \multiput(30,0)(.2,-.112){75}{\circle*{.02}}
\end{picture}
\end{center}
\vspace{1.8cm}
The ideal edge $I=I({\rm cl}(G))$ of the clique clutter of $G$ is not
normal. This graph is chordal, hence perfect. 
Thus edge ideals of clique clutters of perfect graphs are in general
not normal.
\end{example}

\section{TDI systems}\label{tdi-systems}

As already noted in the introduction, the MFMC property for a clutter 
$\CC$ can be stated as $x[A|I_n] \leq (-{\bf 1}|0)$ is TDI, 
where $A$ is the clutter matrix of $\CC$. Thus we have already seen 
(a small slice of) the central role that TDI plays in combinatorial optimization. 
Another motivation to study TDI systems is, for cost vectors $w \in \N^q$ 
for which the problem 
$
{\rm min}\{\langle y,w\rangle \vert\, 
y\geq 0;\, Ay=a \}
$ 
has a unique solution, the system $xA\leq w$ is TDI
if and only if the toric initial ideal of $A$ with respect to $w$ 
is generated by square-free monomials \cite[Corollary. 8.9]{Stur1}. 

In the previous sections we saw how Ehrhart clutters facilitated 
results on the regularity of edge ideals and the MFMC property 
of clutters. In this short closing section we replace the 0/1 clutter matrix 
with a general integer matrix $A$ whose columns, as before, we denote by 
the set $\{ v_1,\ldots, v_q \}$ and we replace the vector 
${\mathbf 1}$ of all ones of length $q$ with a general vector 
$w = (w_1,\ldots, w_q) \in \Z^q$. Motivated by our results 
for Ehrhart clutters we can ask: is there any significance to 
${\mathcal H}(A,w) := \{ (v_1,w_1),\ldots, (v_q,w_q) \}$ 
being a Hilbert basis, especially with regards to the TDI property? 

\begin{theorem}\label{viwi-tdi} Let $A$ be an integral matrix with 
column vectors $v_1,\ldots,v_q$ and let $w=(w_i)$ be an integral vector. 
If the polyhedron $P=\{x\vert\, xA\leq w\}$ is integral 
and $\mathcal{H}(A,w)=\{(v_i,w_i)\}_{i=1}^q$ is a Hilbert 
basis, then the system $xA\leq w$ is TDI.  
\end{theorem}

\begin{proof} To show that $xA\leq w$ is TDI it suffices to show 
that the second part of the equivalence in Theorem~\ref{thm:tdiHilbert} 
holds for the minimal faces of $P$. Let $F$ be such a minimal face of $P$. 
We may assume that $v_1,\ldots,v_r$ are the columns of $A$ which are active in $F$. 
Then $\langle x,v_i\rangle=w_i$ for $x\in F$ and $1\leq i\leq r$. 
If $\langle y,v_i\rangle<w_i$ for some $y\in F$, then $\langle
x,v_i\rangle<w_i$ for any other $x\in F$. Indeed if $\langle
x,v_i\rangle=w_i$ for some $x\in F$, consider the supporting
hyperplane of $P$ given by $H=\{x\vert\langle
x,v_i\rangle=w_i\}$, then $x\in F\cap H\subsetneq F$ because
$y\in F$ and $y\notin F\cap H$, a contradiction to the
minimality of the face $F$. Thus we may also assume that 
$\langle x,v_i\rangle<w_i$ for $x\in F$ and $i>r$. 
Since $P$ is integral, each face of $P$ contains integral vectors, 
see \cite[Section 16.3]{Schr}. Pick an integral vector $x_0\in F$. 
We can now show that $\mathcal{B}=\{v_1,\ldots,v_r\}$ is a Hilbert basis. 

Let ${a}\in\mathbb{R}_+\mathcal{B}\cap\mathbb{Z}^n$. Then we can write 
$a = \lambda_1v_1+\cdots+\lambda_rv_r$ where $\lambda_i \geq 0$ for each 
$i = 1, \ldots, r$. Thus we have
$$
b 
\, := \, 
\langle{a},x_0\rangle 
\, = \, 
\lambda_1\langle v_1,x_0\rangle+ \cdots +\lambda_r\langle v_r,x_0\rangle 
\, = \, \lambda_1w_1+\cdots+\lambda_rw_r.
$$
In particular $b$ is an integer and we can write 
$
(a,b)=\lambda_1(v_1,w_1)+\cdots+\lambda_r(v_r,w_r).
$ 
By hypothesis $\mathcal{H}(A,w)$ is a Hilbert basis. Therefore we can write
$
(a,b)=\eta_1(v_1,w_1)+\cdots+\eta_q(v_q,w_q)$ where 
$
\eta_i \in\mathbb{N}$ 
for each $i = 1, \ldots, q$. 
Therefore
$$
0 \, = \, 
\langle(a,b),(x_0,-1)\rangle
\, = \, 
\sum_{i=1}^r\eta_i\underbrace{\langle(v_i,w_i),(x_0,-1)\rangle}_{=0}+
\sum_{i=r+1}^q\eta_i\underbrace{\langle(v_i,w_i),(x_0,-1)\rangle}_{<0}.
$$
Hence $\eta_i=0$ for $i>r$ and ${a}=\eta_1v_1+\cdots+\eta_r v_r$. 
Thus $a\in\mathbb{N}\mathcal{B}$, as required. 
\end{proof} 

One motivation for the above theorem is \cite[Theorem~22.18]{Schr}, which states 
that the system $xA\leq w$ has the {\em integer rounding property} if and only if the set 
$
\{(v_1,w_1),\ldots,(v_q,w_q),e_{n+1}\}
\subset\mathbb{Z}^{n+1}
$
is a Hilbert basis. A straightforward corollary of this result is that we have an 
``if and only if'' in Theorem~\ref{viwi-tdi} if we append the vector $e_{n+1}$ onto 
${\mathcal H}(A,w)$. The converse of Theorem~\ref{viwi-tdi} is not true in general but 
there are interesting systems for which the converse does hold. Let $A$ be an integral 
matrix and let $w$ be an integral vector. As before, we can rewrite 
the system $x\geq 0;xA\leq w$ as $x[A|-I_n] \leq (w|0)$.

\begin{proposition}\label{tdi-non-negative} 
Let $A$ be a non-negative integral matrix of order
$n\times q$ with column vectors
$v_1,\ldots,v_q$ and let $w=(w_i)\in\mathbb{N}^q$. Then the
system $x\geq 0; xA\leq w$ is TDI if and only if the 
polyhedron $P=\{x\vert\, x\geq 0; xA\leq w\}$ is integral 
and 
${\mathcal H} \, := \, \mathcal{H}([A|-I_n],(w|0)) \, = \, 
\{(v_1,w_1),\ldots,(v_q,w_q),-e_1,\ldots,-e_n\}$ 
is a Hilbert basis. 
\end{proposition}

\begin{proof} The ``if'' part follows from Theorem~\ref{viwi-tdi}. 
For the ``only if'' part, assume that the system $x\geq 0; xA\leq w$ is TDI. 
By the Edmonds-Giles theorem \cite[Corollary~22.1c]{Schr}, the polyhedrom $P$ 
must be integral. All that remains to show is that $\mathcal{H}$ is a Hilbert basis. 
Take $(a,b) \in \mathbb{R}_+{\mathcal H} \cap \mathbb{Z}^{n+1}$, where 
$a \in \mathbb{Z}^n$ and $b\in\mathbb{Z}$. By hypothesis, the linear program 
$
{\rm min}\{\langle y,w\rangle \vert\, y\geq 0; Ay\geq{a}\} 
$
has an integral optimum solution $y=(y_i)$ such that $\langle y,w \rangle \leq b$. 
Since $y\geq 0$ and ${a}\leq Ay$, we can write 
$a \, = \, y_1v_1 + \cdots+y_qv_q-\delta_1e_1-\cdots-\delta_ne_n$ 
where each $\delta_i\in\mathbb{N}$. This implies that 
$$
(a,b)=y_1(v_1,w_1)+\cdots+y_{q-1}(v_{q-1},w_{q-1})+(y_q+b-\langle{y},w\rangle)(v_q,w_q)
-(b-\langle{y},w\rangle)v_q-\delta,
$$
where $\delta=(\delta_i)$. As the entries of $A$ are in $\mathbb{N}$,
the vector $-v_q$ can be written as a non-negative integer
combination of $-e_1,\ldots,-e_n$. Thus
$(a,b)\in\mathbb{N}{\mathcal H}$ as claimed.
\end{proof}

\bibliographystyle{plain}

\end{document}